\documentclass[preprint,12pt]{elsarticle}

\usepackage{amssymb}
\usepackage{amsmath}
\usepackage{amsthm}
\usepackage[colorlinks]{hyperref}
\AtBeginDocument{%
  \hypersetup{%
    linkcolor=blue,%
    citecolor=green,%
  }%
}
\usepackage{cleveref}

\crefname{equation}{}{}

\newtheorem{theorem}{Theorem}[section]
\newtheorem{lemma}[theorem]{Lemma}

\theoremstyle{definition}
\newtheorem{definition}[theorem]{Definition}

\theoremstyle{remark}
\newtheorem{remark}[theorem]{Remark}

\numberwithin{equation}{section}

\journal{}

\begin{document}

\begin{frontmatter}



\title{On the boundary H\"{o}lder regularity for the infinity Laplace equation \tnoteref{t1}}

\author[rvt]{Leyun Wu}
\ead{leyunwu@mail.nwpu.edu.cn}
\author[rvt]{Yuanyuan Lian\corref{cor1}}
\ead{lianyuanyuan@nwpu.edu.cn; lianyuanyuan.hthk@gmail.com}
\author[rvt]{Kai Zhang\corref{cor2}}
\ead{zhang\_kai@nwpu.edu.cn; zhangkaizfz@gmail.com}
\tnotetext[t1]{This research is supported by the National Natural Science Foundation of China (Grant No. 11701454), the Natural Science Basic Research Plan in Shaanxi Province of China (Program No. 2018JQ1039) and the Fundamental Research Funds for the Central Universities (Grant No. 31020170QD032).}

\cortext[cor1]{Corresponding author.}
\cortext[cor2]{ORCID: \href{https://orcid.org/0000-0002-1896-3206}{0000-0002-1896-3206};
MR Author ID: \href{http://mathscinet.ams.org/mathscinet/search/author.html?mrauthid=1098004}
{1098004}.
}

\address[rvt]{Department of Applied Mathematics, Northwestern Polytechnical University, Xi'an, Shaanxi, 710129, PR China}

\begin{abstract}
In this note, we prove the boundary H\"{o}lder regularity for the infinity Laplace equation under a proper geometric condition. This geometric condition is quite general, and the exterior cone condition, the Reifenberg flat domains, and the corkscrew domains (including the non-tangentially accessible domains) are special cases. The key idea, following \cite{ZLH}, is that the strong maximum principle and the scaling invariance imply the boundary H\"{o}lder regularity.
\end{abstract}


\end{frontmatter}

In this note, we prove the boundary H\"{o}lder regularity for the infinity Laplace equation:
\begin{equation} \label{eq1.1}
\left\{\begin{aligned}
     \Delta_\infty u =u_iu_ju_{ij}=0 &\quad \text{~in } \Omega; \\
     u= g & \quad\text{~on } \partial\Omega,
\end{aligned}\right.
\end{equation}
where $\Omega \subset \mathbb{R}^n$ is a bounded domain, $u_i=\partial u/\partial x_i$, $u_{ij}=\partial ^2 u/\partial x_i\partial x_j$ and the Einstein summation convention is used.

The following boundary regularity are well known. If $g\in C^{0,1}(\partial \Omega)$ (i.e., $g$ is Lipschitz continuous on $\partial \Omega$), $u\in C^{0,1}(\bar{\Omega})$ (i.e., $u$ is Lipschitz continuous on $\bar{\Omega}$). In addition, if $g\in C^0(\partial \Omega)$, $u\in C^0(\bar{\Omega})$. It should be noted that both results hold without any geometric condition on $\partial \Omega$.

It is natural to study the boundary H\"{o}lder regularity, i.e., whether $g\in C^{\alpha}(\partial \Omega)$ implies $u\in C^{\alpha}(\bar{\Omega})$. This is the aim of this note. To the authors' knowledge, our result is the first one contributing to the boundary H\"{o}lder regularity.

In this note, we obtain the boundary H\"{o}lder regularity for the solutions of\cref{eq1.1} under a proper geometric condition on $\partial \Omega$. The idea and the method originate from \cite{ZLH}. As pointed out in \cite{ZLH}, this geometric condition is quite general and the exterior cone condition, the Reifenberg flat domains, and the corkscrew domains (including the non-tangentially accessible domains) are special cases. The main idea is that the strong maximum principle implies a decay for the solution, then a scaling argument leads to the H\"{o}lder regularity.

The following is the geometric condition under which we prove the boundary H\"{o}lder regularity.
\begin{definition}\label{Con2.3}\textbf{(Uniform condition)}
Let $\Omega \subset \mathbb{R}^n$ be a bounded domain and $x_0 \in \partial\Omega.$ We say that $\Omega$ satisfies the uniform condition at $x_0$ if the following holds: there exist constants $0<\tau_1<\tau_2<1$, $0<\nu<1$ and a positive sequence $\{r_k\}_{k=0}^{\infty}$ such that
\begin{equation}\label{e-rk}
    \tau_1 r_{k-1}\leq r_k \leq \tau_2 r_{k-1}, ~~\forall ~ k\geq 1,
\end{equation}
and for any $k\geq 0$, there exists $y_k\in \partial B(x_0,r_k)$ such that
\begin{equation}\label{e-bnu}
  \partial B(x_0,r_k) \cap B(y_k,\nu r_k) \subset \Omega ^c.
\end{equation}
\end{definition}

%
%
We say that $g$ is $C^{\alpha}$ at $x_0\in \bar{\Omega}$ or $g\in C^{\alpha}(x_0)$ if there exists a constant $K$ such that
\begin{equation*}
  |g(x)-g(x_0)|\leq K|x-x_0|^{\alpha}, ~\forall ~x\in \bar{\Omega}.
\end{equation*}
Then denote $[g]_{C^{\alpha}(x_0)}=\inf K$.

Our main result is the following.

\begin{theorem}\label{Th1.1}
Suppose that $\Omega$ satisfies the uniform condition at $0\in\partial\Omega$ with $r_0=1$. Let $u$ be a viscosity solution of
\begin{equation*}
\left\{\begin{aligned}
     \Delta_\infty u =0 &\quad \text{~in } \Omega\cap B_1; \\
     u= g & \quad\text{~on } \partial\Omega \cap B_1,
\end{aligned}\right.
\end{equation*}
where $g \in C^\alpha (0).$

Then $u$ is $C^\beta$ at $0$ and
\begin{equation}\label{eq1.2}
  |u(x)-u(0)|\leq 8 |x|^\beta \left( \| u\|_{L^\infty (\Omega \cap B_1)}+[g]_{C^\alpha (0)}\right), \forall x \in \Omega \cap B_1,
\end{equation}
where $0<\beta \leq \alpha$ depends only on $n, \tau_1, \tau_2$ and $\nu$.
\end{theorem}


The proof of Theorem \ref{Th1.1} depends on the solvability and the strong maximum principle for the infinity Laplace equation, which are already known. The following lemma shows the solvability (see Theorem 5.21 in \cite{W-1}).
\begin{lemma}\label{Le2-1}
For any $g \in C^{0,1} (\partial \Omega),$ there exists a unique viscosity solution $u \in C^{0,1} (\bar{\Omega})$ to
 \begin{equation*}  \left\{\begin{aligned}
     \Delta_\infty u =0 &\quad \text{~in } \Omega; \\
     u= g & \quad\text{~on } \partial \Omega.
  \end{aligned}\right.\end{equation*}
  \end{lemma}

The strong maximum principle can be derived easily from the following Harnack inequality (see Theorem 2.21 in \cite{W-1} or Proposition 6.3 in \cite{L}).
\begin{lemma}\label{Le2-2}
Let $u\geq 0$ be a viscosity solution of

\begin{equation*}
\Delta_\infty u =0 \quad \text{~in } B_1.
\end{equation*}
Then
\begin{equation*}
  \mathop {\sup }\limits_{{B_{1/2}}} u \leq 3 \mathop {\inf }\limits_{{B_{1/2}}} u.
\end{equation*}
\end{lemma}

\begin{remark}\label{Re2-3}
Obviously, the Harnack inequality implies the strong maximum principle. Hence, for the solution $u$ of\cref{eq1.1}, we have
\begin{equation*}
  \sup_{\Omega'} u\leq (1-\mu) \sup_{\Omega} u,
\end{equation*}
where $\Omega'\subset\subset \Omega$ and $\mu$ depends only on $\Omega,
\Omega'$ and $g$.
\end{remark}

Let $0<\nu<1$ be as in the uniform condition. Then we choose and fix a function $g_{\nu}\in C^{\infty}(\partial B_1)$ with $0\leq g_{\nu} \leq 1$ and
\begin{equation}\label{gnu}
g_{\nu}(x)\equiv
\left\{ \begin{aligned}
&0 ~~&&\mbox{on}~~\partial B_1\cap B(e_1,\nu/2);\\
&1 ~~&&\mbox{on}~~\partial B_1\backslash B(e_1,\nu),
\end{aligned}\right.
\end{equation}
where $e_1=(1,0,0,...,0)$. Next, introduce functions $g_k$ ($k\geq 0$) with
\begin{equation}\label{gk1}
g_{k}(x)\equiv
\left\{ \begin{aligned}
&0 ~~&&\mbox{on}~~\partial B_{r_k}\cap B(y_k,\nu r_k/2);\\
&1 ~~&&\mbox{on}~~\partial B_{r_k}\backslash B(y_k,\nu r_k)
\end{aligned}\right.
\end{equation}
and
\begin{equation}\label{gk2}
g_k(r_k\cdot T_k x)\equiv g_{\nu}(x)~~\mbox{on}~~\partial  B_1
\end{equation}
for some orthogonal matrix $T_k$. Here $\{r_k\}$ and $\{y_k\}$ are as in the uniform condition. From the strong maximum principle, we have the following simple result.
\begin{lemma}\label{le1}
Let $v$ be a viscosity solution of
\begin{equation*}  \left\{\begin{aligned}
     \Delta_\infty v &=0&&\quad \text{~in } B_{r_k}; \\
     v&= ag_k+b &&\quad\text{~on } \partial B_{r_k},
  \end{aligned}\right.\end{equation*}
where $a,b>0$.

Then
\begin{equation*}\label{eq1.3}
  \sup_{B_{r_{k+1}}}v\leq (1-\mu)a+b,
\end{equation*}
where $0<\mu<1$ depends only on $n,\tau_2$ and $\nu$.
\end{lemma}
\proof Let $y=T^T_kx/r_k$ such that $g_{\nu}(y)=g_k(x)$. Let
\begin{equation*}
  w(y)=\frac{v(x)-b}{a}.
\end{equation*}
Then $w$ satisfies
\begin{equation*}  \left\{\begin{aligned}
     \Delta_\infty w &=0&&\quad \text{~in } B_{1}; \\
     w&= g_{\nu} &&\quad\text{~on } \partial B_{1}.
  \end{aligned}\right.\end{equation*}
Then by the strong maximum principle (see \Cref{Le2-2})
\begin{equation*}
\sup_{B_{\tau_2}}w\leq (1-\mu),
\end{equation*}
where $0<\mu<1$ depends only on $n,\tau_2$ and $g_{\nu}$. Note that $g_{\nu}$ depends only on $\nu$. Hence, $\mu$ depends only on $n,\tau_2$ and $\nu$. By rescaling,
\begin{equation*}
    \sup_{B_{r_{k+1}}}v\leq \sup_{B_{\tau_2 r_{k}}}v\leq (1-\mu)a+b.
\end{equation*}
~\qed

Now we give the

\noindent\textbf{Proof of Theorem \ref{Th1.1}.} We assume that $g(0)=0$. Otherwise, we may consider $v=u-g(0)$. Let $M=\|u\|_{{L^\infty} (\Omega \cap B_1)}+[g]_{C^\alpha(0)}$ and $\Omega_r=\Omega \cap B_r.$ To prove\cref{eq1.2}, we only need to prove the following:

There exists a constant $0<\beta \leq \alpha$ depending only on $n, \tau_1, \tau_2$ and $\nu$ such that
\begin{equation}\label{eq2-1}
\tau_1^\beta \geq \frac{1}{2}
\end{equation}
and for all $k\geq 0,$
\begin{equation}\label{eq2-2}
\|u \|_{L^\infty (\Omega_{r_k})}\leq 4Mr_k^\beta.
\end{equation}

Indeed, suppose that\cref{eq2-1} and\cref{eq2-2} hold. Then for any $x\in \Omega\cap B_{1}$, there exists $k$ such that $r_{k+1}\leq |x| \leq r_{k}$. Hence,
\begin{equation*}
|u(x)|\leq 4Mr_k^{\beta}\leq \frac{4Mr_{k+1}^{\beta}}{\tau_1^{\beta}}\leq 8M|x|^{\beta}.
\end{equation*}

We prove\cref{eq2-2} by induction. For $k=0,$ it holds clearly. Suppose that it holds for $k,$ then we need to prove that it holds for $k+1.$

By \Cref{Le2-1}, there exists a unique viscosity solution $v$ of
 \begin{equation*}  \left\{\begin{aligned}
     \Delta_\infty v =0 &\quad \text{~in } B_{r_k}; \\
     v= \tilde{g} & \quad\text{~on } \partial B_{r_k},
  \end{aligned}\right.\end{equation*}
where $\tilde{g}=(4Mr_k^\beta -Mr_k^\alpha)g_k+Mr_k^\alpha.$  Then it is easy to check that
\begin{equation*}
  -v\leq u \leq v ~\mbox{on}~\partial \Omega_{r_k}.
\end{equation*}
Then by the comparison principle (see \cite[Theorem 5.22]{W-1}),
\begin{equation}\label{eq2-3}
-v\leq u \leq v \text{~in~} \Omega_{r_k}.
\end{equation}
From \Cref{le1}, we have
\begin{equation*}\label{eq2-4}
\begin{aligned}
\sup_{B_{r_{k+1}}}v & \leq (1-\mu)\left(4Mr_k^\beta-Mr_k^\alpha\right)+Mr_k^\alpha  \\
 & = (1-\mu) 4Mr_{k}^\beta+\mu Mr_k^\beta \\
 & \leq 4Mr_{k+1}^\beta \left(\frac{1-\mu}{\tau_1^\beta}+\frac{\mu}{4\tau_1^\beta}\right),
\end{aligned}
\end{equation*}
where $0<\mu<1$ depends only on $n, \tau_2$ and $\nu.$

Take $\beta$ small enough such that\cref{eq2-1} holds and
\begin{equation*}
  \frac{1-\mu}{\tau_1^\beta} <1-\frac{\mu}{2}.
\end{equation*}
Hence,
\begin{equation*}
  \sup_{B_{r_{k+1}}}v\leq 4Mr_{k+1}^\beta.
\end{equation*}
Then, combining with\cref{eq2-3}, we have
\begin{equation*}
  \| u \|_{{L^\infty}(\Omega_{r_{k+1}})}\leq 4 Mr_{k+1}^\beta.
\end{equation*}
By induction, the proof is completed.~\qed

\section*{References}



\end{document}